\documentclass[fleqn,preprint,3p,a4paper]{elsarticle}

\usepackage{amssymb}
\usepackage{amsmath}
\usepackage{amsthm}
\usepackage{dcolumn}
\usepackage{endnotes}
\usepackage{tabularx}
\usepackage[matrix,arrow]{xy}
\usepackage{wasysym}

\newtheorem{theorem}{Theorem}[section]
\newtheorem{proposition}[theorem]{Proposition}
\newtheorem{lemma}[theorem]{Lemma}

\theoremstyle{definition}
\newtheorem{definition}[theorem]{Definition}
\newtheorem{example}[theorem]{Example}
\newtheorem{remark}[theorem]{Remark}

\newcommand{\II}{\begin{enumerate}}
\newcommand{\III}{\end{enumerate}}

\newcommand{\ir}{{\mathsf{Irr}}}

\newcommand{\cl}{{\rm cl}}
\newcommand{\ii}{{\rm int}}

\newcommand{\ua}{\mathord{\uparrow}}
\newcommand{\da}{\mathord{\downarrow}}

\newcommand{\wdd}{\mathord{\mathsf{WD}}}

\journal{}

\begin{document}

\begin{frontmatter}



\title{On three problems about well-filteredness of $T_0$-spaces\tnoteref{t1}}
\tnotetext[t1]{This research was supported by the National Natural Science Foundation of China (Nos. 12471070, 12071199)}

\author[X. Xu]{Xiaoquan Xu\corref{mycorrespondingauthor}}
\cortext[mycorrespondingauthor]{Corresponding author}
\ead{xiqxu2002@163.com}
\address[X. Xu]{School of Mathematics and Statistics, Minnan Normal University, Zhangzhou 363000, China}

\author[X. Wen]{Xinpeng Wen}
\ead{wenxinpeng2009@163.com}
\address[X. Wen]{College of Mathematics and Information, Nanchang Hangkong University, Jiangxi 330063, China}

\begin{abstract}
In this paper, we show that there is a countable Noetherian complete lattice $L$ and an order-compatible $d$-topology $\tau$ on $L$ such that $(L, \tau)$ is not well-filtered, and there exist a dcpo $P$ and an order-compatible well-filtered topology $\tau$ on $P$ but the Scott topology $\sigma (P)$ is not well-filtered. For such poset $P$ and topology $\tau$, let $Y=(P, \tau)$ and $X = 1$ (the topological
space with single point), then the function space $\mathbb{C}(X, Y)$ equipped with the Scott topology is not well-filtered. These results  answer three open problems concerning the well-filteredness of $T_0$-spaces.
\end{abstract}


\end{frontmatter}

In domain theory, the $d$-spaces and sober spaces are two important classes of $T_0$-spaces (see \cite{Ershov-1999, GHKLMS-2003, Goubault-2013, Wyler-1981, Xu-Zhao-2021}). In \cite{Heckmann-1992}, Heckmann introduced a new types $T_0$-spaces ---  well-filtered spaces, which lie between $d$-spaces and sober spaces. The well-filteredness has emerged as a very useful and important property for $T_0$-spaces (see [5, 7, 11, 12, 14-20, 22, 24-38]).

In this paper, we give answers to three open problems concerning the well-filteredness of $T_0$-spaces. The paper is organized as follows.

Section 1 provides some fundamental definitions and notations about topology and lattice-ordered structures which will be used in the whole paper. Also the filtered version of topological Rudin lemma and some characterizations of $d$-spaces and well-filtered spaces are listed.

Section 2 gives a countable Noetherian complete lattice $L$ and an order-compatible $d$-topology $\tau$ on $L$ such that $(L, \tau)$ is not well-filtered, which gives an answer to \cite[Question 4.7]{Xu-2026-2}.

In Section 3, we show that if $\mathbb{J}_{\top}$ is the dcpo obtaining by adding a top element $\top$ to Johnstone's dcpo $\mathbb{J}$ and $\tau=\{\emptyset\}\cup \{U\cup\{\top\} : U\in \sigma(\mathbb{J})\setminus \{\emptyset\}\}$, then  $(\mathbb{J}_{\top}, \tau)$ is sober (and hence well-filtered) but the Scott space $\Sigma~\!\!\mathbb{J}_\top$ is not well-filtered, which answers \cite[Question 64 and Question 67]{Xu-Bao-Zhang-2022}.  Furthermore, we show that the function space $\mathbb{C}(1, (\mathbb{J}_{\top}, \tau))$ ($1$ is the topological space with single point) equipped with the Scott topology is not well-filtered, which gives an answer to \cite[Question 70]{Xu-Bao-Zhang-2022}.

\section{Preliminaries}

For a poset $P$ and $A\subseteq P$, let
$\mathord{\downarrow}A=\{x\in P: x\leq  a \mbox{ for some }
a\in A\}$ and $\mathord{\uparrow}A=\{x\in P: x\geq  a \mbox{
	for some } a\in A\}$. For  $x\in P$, we write
$\mathord{\downarrow}x$ for $\mathord{\downarrow}\{x\}$ and
$\mathord{\uparrow}x$ for $\mathord{\uparrow}\{x\}$.  A subset $A$
is called a \emph{lower set} (resp., an \emph{upper set}) if
$A=\mathord{\downarrow}A$ (resp., $A=\mathord{\uparrow}A$). Let $P^{(<\omega)}=\{F\subseteq P : F \mbox{~is a nonempty finite set}\}$ and $\mathbf{Fin}~\!P=\{{\uparrow} F : F\in P^{(<\omega)}\}$. 
An element $u$ is said to be an \emph{upper bound} of a set $C$, if $c\leq u$ for all $c\in C$ (i.e., $C\subseteq {\downarrow} u$). For $E\subseteq P$, if the set of upper bounds of $E$ in $P$ has a unique smallest element, we call this element the \emph{least upper bound} and write it as $\vee E$ or sup~\!\! $E$ (for \emph{supremum}). Dually, the \emph{greatest lower bound} of $E$ in $P$ is written as $\wedge E$ or inf~\!\! $E$ (for \emph{infimum}). The poset $P$ is called an \emph{inf semilattice} (shortly \emph{semilattice}) if for any two elements
$a, b\in P$, $a\wedge b$ exists in $P$. Dually, $P$ is a \emph{sup semilattice} if for any two elements $a, b\in P$, $a\vee b$ exists in $P$.
The set of all natural
numbers is denoted by $\mathbb{N}$. Let $\mathbb{N}^+=\mathbb{N}\setminus \{0\}$. For a set $X$, let $|X|$ be the cardinality of $X$ and $2^X$ the set of all subsets of $X$.

A nonempty subset $D$ of a poset $Q$ is \emph{directed} if every two
elements in $D$ have an upper bound in $D$. The set of all directed sets of $Q$ is denoted by $\mathcal D(Q)$.  The poset $Q$ is called a \emph{directed complete poset}, or \emph{dcpo} for short, if for any
$D\in \mathcal D(Q)$, its supremum $\vee D$ exists in $Q$. The notion $x=\bigvee^{\uparrow}G$ is a convenient device to express that, firstly, the set $G$ is directed and, secondly, $x$ is its leat upper bounded.

A subset $U$ of a poset $P$ is \emph{Scott open} if
(i) $U=\mathord{\uparrow}U$, and (ii) for any directed subset $D$ for
which $\vee D$ exists, $\vee D\in U$ implies $D\cap
U\neq\emptyset$. All Scott open subsets of $P$ form a topology,
and we call this topology  the \emph{Scott topology} on $P$ and
denote it by $\sigma(P)$. The space $\Sigma~\!\! P=(P,\sigma(P))$ is called the
\emph{Scott space} of $P$. The \emph{upper topology} on a poset $P$, generated
by $\{P\setminus {\downarrow}x : x\in P\}$ (as a subsase), is denoted by $\upsilon (P)$. The upper sets
form the (\emph{upper}) \emph{Alexandroff topology} $\alpha (P)$.

For a $T_0$-space $X$, we use $\mathcal O(X)$ (resp., $\Gamma(X)$) to denote the set of all open subsets (resp., closed subsets) of $X$. The closure of a subset $A$ in $X$ is denoted by $\cl_X A$ (or simply by $\cl A$ if there is no ambiguity) or $\overline{A}$, and the interior of $A$ will be denoted by $\ii_X A$ or simply by $\ii A$. We use $\leq_X$ to represent the \emph{specialization order} of $X$, that is, $x\leq_X y$ if{}f $x\in \overline{\{y\}}$.  We will use $\Omega~\!\!X$ or even $X$ to denote the poset $(X, \leq_X)$. In what follows, when a $T_0$-space is considered as a poset, the order always refers to the specialization order if no other explanation is given. Let $\mathcal S_c(X)=\{\overline{{\{x\}}} : x\in X\}$ and $\mathcal D_c(X)=\{\overline{D} : D\in \mathcal D(X)\}$. For a poset $P$, a $T_0$-topology $\tau$ on $P$ is said to be \emph{order}-\emph{compatible} if $\leq_{\tau}$ agrees with the original order on $P$. It is easy to verify that $\tau$ is order-compatible iff $\upsilon (P)\subseteq \tau\subseteq \alpha (P)$. 

Let $\mathbf{Set}$ denote the category of all sets and mappings. The category of all $T_0$-spaces with continuous mappings is denoted by $\mathbf{Top}_0$. A $T_0$-space $X$ is called a \emph{$d$-space} (or \emph{monotone convergence space}) if $X$ (with the specialization order) is a dcpo and $\mathcal O(X) \subseteq \sigma(X)$ (cf. \cite{GHKLMS-2003, Wyler-1981}). For a set $Y$ and a topology $\tau$ on $Y$, we call $\tau$ a $d$-\emph{topology} on $Y$ if $(Y, \tau)$ is a $d$-space. A topology $\delta$ on a poset $P$ is said to be an \emph{order}-\emph{compatible} $d$-\emph{topology} if $\delta$ is order-compatible and $(P, \delta)$ is a $d$-space. Clearly, $\delta$ is an  order-compatible $d$-topology on $P$ iff $P$ is a dcpo and $\upsilon (P)\subseteq \delta \subseteq \sigma (P)$.

\begin{proposition}\label{prop-d-space-charac} (\cite[Proposition 3.3]{Xu-Shen-Xi-Zhao-2020-2}) For a $T_0$-space $X$, the following two conditions are equivalent:
\begin{enumerate}[\rm (1)]
\item $X$ is a $d$-space.
           \item $\mathcal D_c(X)=\mathcal S_c(X)$.
\end{enumerate}
\end{proposition}

A nonempty subset $A$ of a $T_0$-space $X$ is said to be \emph{irreducible} if for any $\{F_1, F_2\}\subseteq \Gamma(X)$, $A \subseteq F_1\cup F_2$ implies $A \subseteq F_1$ or $A \subseteq  F_2$.  Denote by $\ir(X)$ (resp., $\ir_c(X)$) the set of all irreducible (resp., irreducible closed) subsets of $X$. Clearly, every directed subset of $X$ (with the specialization order) is irreducible. The space $X$ is called \emph{sober}, if for any  $A\in\ir_c(X)$, there is a unique point $x\in X$ such that $A=\overline{\{x\}}$.

A subset $A$ of a $T_0$-space $X$ is called \emph{saturated} if it equals the intersection of all open sets containing it (equivalently, $A$ is an upper set with respect to the specialization order). We use $\mathord{Q}(X)$ to
denote the set of all nonempty compact saturated subsets of $X$.

\begin{definition}\label{def-WF-space} A $T_0$-space $X$ is said to be \emph{well-filtered} if for any open set $U$ and filtered family $\mathcal{K}\subseteq \mathord{Q}(X)$, $\bigcap\mathcal{K}{\subseteq} U$ implies $K{\subseteq} U$ for some $K{\in}\mathcal{K}$.
\end{definition}

The well-filtered space was originally introduced and investigated by Heckmann \cite{Heckmann-1992}, which was called the $\mathcal U_K$-\emph{admitting space} in \cite{Heckmann-1992}.

Rudin's Lemma, due to Mary Rudin \cite{Rudin-1981}, is an important tool in domain theory and non-Hausdorff topology (see \cite{Abramsky-Jung-1994, GHKLMS-2003, Gierz-Lawson-Stralka-1983, Goubault-2013}). In \cite{Heckmann-Keimel-2013}, Heckmann and Keimel presented a topological variant of Rudin's Lemma.

In this paper we only need the following filtered version of topological Rudin lemma.

\begin{lemma}\label{t Rudin} Let $X$ be a $T_0$-space and $\mathcal K$ a filtered family of nonempty compact saturated sets of $X$. Then every closed set $C$ of $X$  that
meets all members of $\mathcal{K}$ contains an minimal irreducible closed subset $A$ that meets all
members of $\mathcal{K}$.
\end{lemma}

For a $T_0$-space $X$ and $\mathcal{K}\subseteq \mathord{Q}(X)$, let $M(\mathcal{K})=\{A\in \Gamma(X) : K\cap A\neq\emptyset \mbox{~for all~} K\in \mathcal{K}\}$ and $m(\mathcal{K})=\{A\in \Gamma(X) : A \mbox{~is a minimal menber of~} M(\mathcal{K})\}$.

\begin{definition}\label{rudinset} (\cite[Definition 2.1]{Shen-Xi-Xu-Zhao-2019} and \cite[Definition 4.6 and Definition 6.1]{Xu-Shen-Xi-Zhao-2020-2})
		Let $X$ be a $T_0$-space and $A$ a nonempty closed set of $X$.
\begin{enumerate}[\rm (1)]
\item $A$ is said to be  a \emph{Rudin set}, if there exists a filtered family $\mathcal K\subseteq \mathord{Q}(X)$ such that $\overline{A}\in m(\mathcal K)$ (that is, $\overline{A}$ is a minimal closed set that intersects all members of $\mathcal K$). Denote by $\mathsf{RD}(X)$ (resp., $\mathsf{RD}_c(X)$) the set of all Rudin sets (resp., all closed Rudin sets) of $X$. Clearly, $\mathsf{RD}_c(X)=\{\overline{A} : A\in \mathsf{RD}(X)\}$.
\item $A$ is said to be  a \emph{well-filtered determined set}, $\wdd$ \emph{set} for short, if for any continuous mapping $ f:X\longrightarrow Y$
into a well-filtered space $Y$, there exists a unique $y_A\in Y$ such that $\overline{f(A)}=\overline{\{y_A\}}$.
Denote by $\mathsf{WD}(X)$ (resp., $\mathsf{WD}_c(X)$) the set of all well-filtered determined subsets (resp., all closed well-filtered determined subsets) of $X$. It is easily verified that $\mathsf{WD}_c(X)=\{\overline{A} : A\in \mathsf{WD}(X)\}$.
\end{enumerate}
\end{definition}

\begin{lemma}\label{DRWIsetrelation} (\cite[Proposition 6.2]{Xu-Shen-Xi-Zhao-2020-2})
	Let $X$ be a $T_0$-space. Then $S_c(X)\subseteq\mathcal{D}_c(X)\subseteq \mathsf{RD}(X)\subseteq\mathsf{WD}(X)\subseteq\ir_c(X)$.
\end{lemma}

Using Rudin sets and WD-sets, we have the following characterization of well-filtered spaces.

\begin{proposition}\label{prop-well-filtered-charac} \emph{(\cite[Corollary 7.11]{Xu-Shen-Xi-Zhao-2020-2})} Let $X$ be a $T_{0}$-space. Then the following conditions are equivalent:
\begin{enumerate}[\rm (1)]
\item $X$ is well-filtered.
\item $\mathsf{WD}_{c}(X)=\mathcal{S}_{c}(X)$.
\item $\mathsf{RD}_{c}(X)=\mathcal{S}_{c}(X)$.
\end{enumerate}
\end{proposition}

By Proposition \ref{prop-d-space-charac}, Lemma \ref{DRWIsetrelation} and Proposition \ref{prop-well-filtered-charac}, we have the following:

\begin{proposition}\label{prop-sober-WF-d-space} Every sober space is well-filtered and every well-filtered space is a $d$-space.
\end{proposition}

A poset $P$ is said to be \emph{Noetherian} if it satisfies the \emph{ascending chain condition} ($\mathrm{ACC}$ for short): every ascending chain has a greatest member.

\begin{proposition}\label{prop-Noethrian-Alexandroff-topology} (\cite[Proposition 5.4 and Theorem 5.7]{Zhao-Ho-2015} and \cite[Proposition 3.8]{Xu-Shen-Xi-Zhao-2020-1}) For a poset $P$, the following conditions are equivalent:
	\begin{enumerate}[\rm (1)]
		\item $P$ is Noetherian.
\item Every directed subset of $P$ has a largest member.
        \item $P$ is a dcpo and ${\uparrow}x\in \sigma(P)$ for all $x\in P$.
        \item $P$ is a dcpo and $\sigma(P)=\alpha(P)$.
        \item $(P,\alpha(P))$ is a $d$-space.
        \item $(P,\alpha(P))$ is well-filtered.
\item $(P,\alpha(P))$ is sober.
	\end{enumerate}
\end{proposition}

In order to provide a uniform approach to $d$-spaces, sober spaces and well-filtered spaces and develop a general framework for dealing with all these spaces, Xu \cite{Xu-2021-1} introduced the following concepts.

\begin{definition}\label{def-R-subset-system} (\cite[Definition 3.1 and Definition 3.2]{Xu-2021-1}) (1) A covariant functor ${\rm H}: \mathbf{Top}_0 \rightarrow \mathbf{Set}$ is called a \emph{subset system} on $\mathbf{Top}_0$ provided that the following two conditions are satisfied:
\begin{enumerate}[\rm (i)]
\item $\mathcal S(X)\subseteq {\rm H}(X)\subseteq 2^X$ (the set of all subsets of $X$) for each $X\in$ \emph{ob}($\mathbf{Top}_0$).
\item For any continuous mapping $f: X \rightarrow Y$ in $\mathbf{Top}_0$, ${\rm H}(f)(A)=f(A)\in {\rm H}(Y)$ for all $A\in {\rm H}(X)$.
\end{enumerate}
(2) A subset system ${\rm H}: \mathbf{Top}_0 \rightarrow \mathbf{Set}$ is called an \emph{irreducible subset system}, or an \emph{R-subset system} for short, if ${\rm H}(X)\subseteq \ir (X)$ for all $X\in$ \emph{ob}($\mathbf{Top}_0$).
\end{definition}

In what follows, the capital letter H always stands for an R-subset system ${\rm H} : \mathbf{Top}_0 \rightarrow \mathbf{Set}$. For a $T_0$-space $X$, let ${\rm H}_c(X)=\{\overline{A} : A\in {\rm H}(X)\}$. 

Here are some important examples of R-subset systems used in this paper:

\begin{enumerate}[\rm (1)]
    \item $\mathcal S$ (for $X\in$ \emph{ob}($\mathbf{Top}_0$), $\mathcal S(X)$ is the set of all single point subsets of $X$).
    \item $\mathcal C$ (for $X\in$ \emph{ob}($\mathbf{Top}_0$), $\mathcal C(X)$ is the set of all chains of $X$).
    \item $\mathcal D$ (for $X\in$ \emph{ob}($\mathbf{Top}_0$), $\mathcal D(X)$ is the set of all directed subsets of $X$).
    \item $\mathsf{RD}$ (for $X\in$ \emph{ob}($\mathbf{Top}_0$), $\mathsf{RD}(X)$ is the set of all Rudin subsets of $X$).
    \item $\mathsf{WD}$ (for $X\in$ \emph{ob}($\mathbf{Top}_0$), $\mathsf{WD}(X)$ is the set of all well-filtered determined subsets of $X$).
    \item $\mathcal R$ (for $X\in$ \emph{ob}($\mathbf{Top}_0$), $\mathcal R(X)$ is the set of all irreducible subsets of $X$).
\end{enumerate}

\begin{definition}\label{def-H-sober} (\cite[Definition 4.1]{Xu-2021-1}) Let ${\rm H} : \mathbf{Top}_0 \rightarrow \mathbf{Set}$ be an R-subset system and $X$ be a $T_{0}$-space. $X$ is called \emph{H}-\emph{sober} if for any $A\in {\rm H}(X)$, there is a (unique) point $x\in X$ such that $\overline{A}=\overline{\{x\}}$ or, equivalently, if ${\rm H}_c(X)=\mathcal S_c(X)$.
\end{definition}

Sober spaces are precisely $\mathcal R$-sober spaces. By Proposition \ref{prop-d-space-charac} and Proposition \ref{prop-well-filtered-charac}, $\mathcal D$-sober spaces and $\mathsf{RD}$-sober spaces (or  $\mathsf{WD}$-sober spaces) are exactly $d$-spaces and well-filtered spaces respectively.

\section{Order-compatible $d$-topologies on complete lattices}

For the sobriety of upper topology on a poset, we have the following result.

\begin{proposition}\label{prop-upper-topology-sober} (\cite[Proposition 2.9]{Xu-Shen-Xi-Zhao-2020-2}) For a poset $P$, the space $(P, \upsilon (P))$ is sober iff  $\vee A$ exists in $P$ for all $A\in\ir ((P, \upsilon (P))$. Therefore, for any complete lattice $L$, $(L, \upsilon (L))$ is sober.
\end{proposition}

In \cite{Xi-Lawson-2017}, the following useful result was given.

\begin{proposition}\label{prop-Scott-topology-on-complete-lattice-WF} (\cite[Corollary 3.2]{Xi-Lawson-2017}) If $L$ is a complete lattice, then $\Sigma~\!\! L$ is well-filtered.
\end{proposition}

By Proposition \ref{prop-sober-WF-d-space}, Proposition \ref{prop-upper-topology-sober} and Proposition \ref{prop-Scott-topology-on-complete-lattice-WF}, the upper topology $\upsilon (L)$ and Scott topology $\sigma (L)$ on a complete lattice $L$ are well-filtered. One might naturally wonder whether any order-compatible $d$-topology on $L$ is well-filtered. So in \cite[Question 4.7]{Xu-2026-2} Xu posed the following question:

 \vspace{2mm}

\noindent {\bf Questian A.} Is every order-compatible $d$-topology on a complete lattice a well-filtered topology? Or equivalently, is there a complete lattice $L$ and a topology $\upsilon (L) \varsubsetneqq \tau \varsubsetneqq \sigma (L)$ such that $(L, \tau)$ is not well-filtered?

In the following example we give a countable Noetherian complete lattice $L$ and an order-compatible $d$-topology $\tau$ on $L$ such that $(L, \tau)$ is not well-filtered, which gives an answer to Question A.

\begin{example}\label{exam-order-compatible-d-topology-on-complete-lattice-not-WF} Let $L=\mathbb{N}^+\cup\{\top, \bot\}$ and define an order on $L$ as follows (see Figure 1):
\begin{enumerate}[\rm (i)]
\item $\bot <n<\top$ for any $n\in\mathbf{N}^+$, and
\item $n$ and $m$ are incomparable for any $n, m\in \mathbf{N}^+$ with $n\neq m$.
\end{enumerate}

\begin{figure}[ht]
	\centering
	\includegraphics[height=4cm,width=6cm]{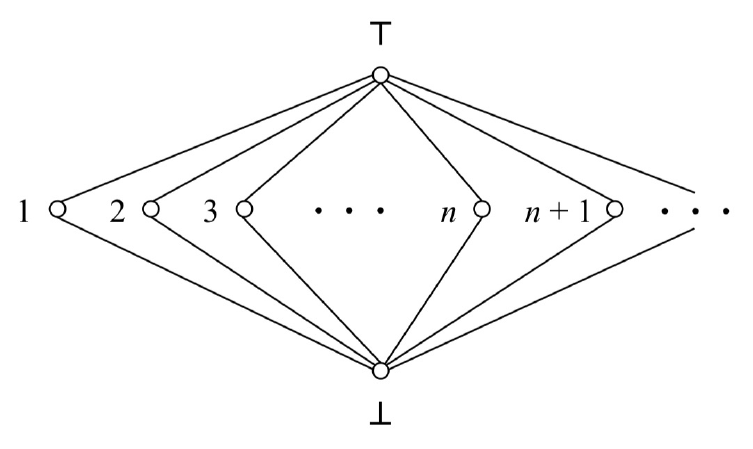}
	\caption{The complete lattice $L$ in Example \ref{exam-order-compatible-d-topology-on-complete-lattice-not-WF}}
\end{figure}

Then $L$ is a countable complete lattice and Noetherian. Let $\tau=\{U\in \alpha (L) : \mathbb{N}^+\cap U\neq\emptyset  \mbox{ implies } \{k\in \mathbb{N}^+ : k\geq m\}\subseteq U \mbox{ for some } m\in\mathbb{N}^+\}$. It is straightforward to verify that $\tau$ is a $T_0$-topology on $L$ and $\upsilon(L) \subseteq \tau$.

Considering the $T_0$-space $(L, \tau)$, we have the following conclusions:

\begin{enumerate}[\rm (a)]
\item $\upsilon(L) \subsetneqq \tau \subsetneqq \sigma(L)$.

As $\{\top\}\in \tau$ and $\{\top\}\notin \upsilon(L)$, we have $\upsilon(L) \subsetneqq \tau$. By Proposition \ref{prop-Noethrian-Alexandroff-topology}, $\tau\subseteq \alpha(L)=\sigma(L)$. Clearly, ${\uparrow}1=\{1, \top\}\in \sigma(L)$, but $\{1, \top\}\notin \tau$. Hence $\tau \subsetneqq \sigma(L)$.

\item $Q((L,\tau))=\alpha (L)\setminus \{\emptyset\}$.

By (a), we have $Q((L,\tau))\subseteq\alpha (L)\setminus \{\emptyset\}$. Conversely, assume that $V$ is a nonempty upper set of $L$. If $V$ is finite, then $V\in Q((L,\tau))$. Now suppose that $V$ is infinite and $\{U_i : i\in I\}\subseteq \tau$ is an open cover of $V$. Then $\emptyset \neq \mathbb{N}^+\cap V\subseteq \mathbb{N}^+\cap \bigcup_{i\in I}U_i=\bigcup_{i\in I}\mathbb{N}^+\cap U_i$. Hence there exists $i_0\in I$ such that $\mathbb{N}^+\cap U_{i_0}\neq\emptyset$, and consequently, $\{k\in \mathbb{N}^+ : k\geq m\}\subseteq U_{i_0}$ for some $m\in \mathbb{N}^+$. As $V\setminus \{k\in \mathbb{N}^+ : k\geq m\}$ is a finite set, there is $I_0\in I^{(<\omega)}$ such that  $V\setminus \{k\in \mathbb{N}^+ : k\geq m\}\subseteq \bigcup_{i\in I_0}U_i$. Then $V\subseteq \bigcup_{i\in I_0\cup\{i_0\}}U_i$. Thus $V\in Q((L,\tau))$, and this completes the proof that $Q((L,\tau))=\alpha (L)\setminus \{\emptyset\}$.

\item~$(L,\tau)$ is not well-filtered.

For $n\in \mathbb{N}^+$, let $A_n=\{l\in \mathbb{N}^+ : l\geq n\}$. Then $\{\ua A_n : n\in \mathbb{N}^+\}\subseteq Q((L, \tau))$ is a filtered family and $\bigcap_{n\in \mathbb{N}^+}\ua A_n=\{\top\}\in \tau$. But for any $m\in \mathbb{N}^+$, $\ua A_m\nsubseteq \{\top\}$. Therefore, $(L,\tau)$ is not well-filtered.
\end{enumerate}
\end{example}

\section{Function spaces equipped with Scott topology}

Function spaces (equipped with certain topologies) are important structures in topology and domain theory (see \cite{Engelking-1989, GHKLMS-2003}). Cartesian closed categories of domains (more generally, certain topological spaces) are appropriate for models of various typed and untyped lambda-calculi and functional programming languages (see \cite{GHKLMS-2003}). Since whether certain properties of topological spaces are preserved when passing to function spaces is connected with the cartesian closed category of topological spaces, this question has attracted considerable attention in domain theory and non-Hausdorff topology, especially for domains (endowed with the Scott topology), sober spaces, $d$-spaces and well-filtered spaces (see \cite{Ershov-Schwidefsky-2020, GHKLMS-2003, Liu-Li-Ho-2021, Xu-2021-1}).

For a topological space $X$ and a $T_0$-space $Y$, let $\mathbb{C}(X, Y)$ denote the set of all continuous functions from $X$ to $Y$. As $Y$ (with the specialization order) is a poset, $\mathbb{C}(X, Y)$ is a poset with the pointwise order ($f\leq g$ iff $f(x)\leq_Y g(x)$ for all $x\in X$). Denote by $[X\rightarrow Y]_{\Sigma}$ the Scott space $\Sigma \mathbb{C}(X, Y)$. In the following, for $y\in Y$, $c_y$ denotes the constant function from $X$ to $Y$ with value $y$ (i.e., $c_{y}(x)=y$ for all $x\in X$). For a point $x\in X$ and an open set $V\in \mathcal{O}(Y)$, let $S(x, V)=\{f\in \mathbb{C}(X, Y)$ : $f(x)\in V\}$. The set $\{S(x, V): x\in X, V\in \mathcal{O}(Y)\}$ is a subbasis for the \emph{pointwise convergence topology} (i.e.,
the relative product topology) on $\mathbb{C}(X, Y)$.

\begin{lemma}\label{lem-function-space-basic-1} (\cite[Lemma 23]{Xu-Bao-Zhang-2022}) Let $X, Y$ be $T_0$-spaces and $\emptyset\neq\mathcal {F}\subseteq \mathbb{C}(X, Y)$. For any mapping $g : X \rightarrow Y$, the following two conditions are equivalent:
\begin{enumerate}[\rm (1)]
\item For each $x\in X$, $\overline{\{f(x) : f\in \mathcal{F}\}}=\overline{\{g(x)\}}$.
\item For each $V\in \mathcal O(X)$, $g^{-1}(V)=\bigcup\limits_{f\in \mathcal F}f^{-1}(V)$.
\end{enumerate}
Therefore, when condition (1) is satisfied, we have that $g\in \mathbb{C}(X, Y)$, $g(x)=\bigvee_{f\in \mathcal{F}} f(x)$ for each $x\in X$ and $g=\bigvee_{\mathbb{C}(X, Y)}\mathcal{F}$.
\end{lemma}

\begin{definition}\label{def-property-S} (\cite[Definition 26]{Xu-Bao-Zhang-2022}) Let ${\rm H} : \mathbf{Top}_0 \rightarrow \mathbf{Set}$ be an R-subset system and $X, Y$ $T_{0}$-spaces. An order-compatible topology $\mathcal T$ on $\mathbb{C}(X, Y)$ is said to have \emph{property S} with respect to H if it satisfies the following condition:

(S)   For any $\mathcal F\in {\rm H}((\mathbb{C}(X, Y), \mathcal T))$ and $g\in \mathbb{C}(X, Y)$, if $g^{-1}(V)=\bigcup_{f\in\mathcal F}f^{-1}(V)$ for each $V\in \mathcal O(Y)$, then $g\in \cl_{\mathcal T}\mathcal F$.
\end{definition}

\begin{remark}\label{rem-property S} By the order-compatibility of $\mathcal T$ and Lemma \ref{lem-function-space-basic-1}, $\mathcal T$ has property S with respect to H iff for any $\mathcal F\in {\rm H}((\mathbb{C}(X, Y), \mathcal T))$ and $g\in \mathbb{C}(X, Y)$ with $g^{-1}(V)=\bigcup_{f\in\mathcal F}f^{-1}(V)$ for all $V\in \mathcal O(Y)$, we have $\cl_{\mathcal T}\mathcal F=\cl_{\mathcal T}\{g\}$.
\end{remark}

\begin{theorem}\label{theor-H-sober-property-S-1} (\cite[Theorem 36]{Xu-Bao-Zhang-2022}) Let ${\rm H} : \mathbf{Top}_0 \rightarrow \mathbf{Set}$ be an R-subset system and $X$ a $T_{0}$-space. Suppose that $\mathcal T$ is a topology on $\mathbb{C}(X, Y)$ which is finer than the pointwise convergence topology. Consider the following two conditions:

\begin{enumerate}[\rm (1)]
\item $(\mathbb{C}(X, Y), \mathcal T)$ is H-sober.
\item $\mathcal T$ has property S with respect to H.
\end{enumerate}
Then (1) $\Rightarrow$ (2). Moreover, if $Y$ is H-sober, then the two conditions are equivalent.
\end{theorem}

As an important corollary of Theorem \ref{theor-H-sober-property-S-1}, we have the following.

\begin{theorem}\label{theor-Scott-functions-H-sober-1} (\cite[Theorem 61]{Xu-Bao-Zhang-2022}) Let ${\rm H} : \mathbf{Top}_0 \rightarrow \mathbf{Set}$ be an R-subset system, $X$ a $T_{0}$-space and $Y$ a $d$-space. Consider the following two conditions:
\begin{enumerate}[\rm (1)]
\item $[X\rightarrow Y]_\Sigma$ is H-sober.
\item The Scott topology on $\mathbb{C}(X, Y)$ has property S with respect to H.
\end{enumerate}
Then (1) $\Rightarrow$ (2). Moreover, if $Y$ is H-sober, then the two conditions are equivalent.
\end{theorem}

\begin{proposition}\label{prop-Scott-is-S-with-respect-to-D}  (\cite[Proposition 62 and Proposition 65]{Xu-Bao-Zhang-2022}) Let $X, Y$ be $T_{0}$-spaces. Then, the Scott topology on $\mathbb{C}(X, Y)$ has property S with respect to $\mathcal D$. Therefore, when $Y$ is a $d$-space, the function space $[X\rightarrow Y]_\Sigma $ is a $d$-space.
\end{proposition}

\begin{proposition}\label{prop-Scott-is-not-S-with-respect-to-R} (\cite[Proposition 63]{Xu-Bao-Zhang-2022}) Let $X$ be any $T_0$-space for which $\Sigma \mathcal O(X)$ is non-sober. Then
\begin{enumerate}[\rm (1)]
\item $[X\rightarrow \Sigma 2]_\Sigma$ is not sober.
\item the Scott topology on $\mathbb{C}(X, \Sigma 2)$ does not have property S with respect to $\mathcal R$.
\end{enumerate}
\end{proposition}

There exist even two sober spaces $X$ and $Y$ such that the Scott topology on $\mathbb{C}(X, Y)$ does not have property S with respect to $\mathcal R$, as shown in the following example.

 \begin{example}\label{exam-Scott-not-S-with-respect-to-R}
 Let $X=1$ be the topological space with single point and $L$ the complete lattice constructed by Isbell in \cite{Isbell-1982}. It is well-known that $\Sigma L$ is non-sober. Let $Y=(L, \upsilon (L))$. Then by Proposition \ref{prop-upper-topology-sober}, $Y$ is sober. Clearly,  $\mathbb{C}(X, Y)=\{c_y : y\in Y\}$ and $c_y\mapsto y : [X\rightarrow Y]_{\Sigma}\rightarrow \Sigma \Omega Y=\Sigma L$ is a homeomorphism. So $Y$ is sober but the function space $[X\rightarrow Y]_\Sigma$ is non-sober, and hence by Theorem \ref{theor-Scott-functions-H-sober-1} (for $H=\mathcal R$), the Scott topology on $\mathbb{C}(X, Y)$ does not have property S with respect to $\mathcal R$.
\end{example}

It is still not known whether a similar result to Proposition \ref{prop-Scott-is-S-with-respect-to-D} holds for $\mathsf{WD}$ (or $\mathsf{RD}$) and well-filtered spaces. So in \cite{Xu-Bao-Zhang-2022} Xu et al. posed the following question (see \cite[Question 64 and Question 67]{Xu-Bao-Zhang-2022}).

\vspace{2mm}

\noindent {\bf Question B.} For $T_0$-spaces $X$ and $Y$, whether the Scott topology on $\mathbb{C}(X, Y)$ has property S with respect to $\mathsf{WD}$ or $\mathsf{RD}$? Especially, for a $T_0$-space $X$ and a well-filtered space $Y$, whether the Scott topology on $\mathbb{C}(X, Y)$ has property S with respect to $\mathsf{WD}$ or $\mathsf{RD}$? Or equivalently, is the function space $[X\rightarrow Y]_\Sigma$ well-filtered?

\vspace{2mm}

The following related question was also posed in \cite[Question 70]{Xu-Bao-Zhang-2022}.

\vspace{2mm}

\noindent {\bf Question C.} Is there a well-filtered space $Y$ such that $\Sigma Y$ (i.e., $\Sigma \Omega Y$) is not well-filtered? Or equivalently, is there a dcpo $P$ and a topology $\upsilon (P)\subseteq \tau\subseteq \sigma (P)$ such that $(P, \tau)$ is well-filtered but $(P, \sigma (P))$ is not well-filtered?

\begin{remark}\label{rem-problem-3-and-problem-2} If the answer of Question C is ``Yes", then the answer of Question B is ``No"! Indeed, suppose that $Y$ is any well-filtered space for which the Scott space $\Sigma Y$ is not well-filtered and $X=1$ (the topological space with single point). Then the function space $\mathbb{C}(X, Y)$ equipped with the Scott topology is not well-filtered since $[X\rightarrow Y]_{\Sigma}$ and $\Sigma Y$ are homeomorphic (see Example \ref{exam-Scott-not-S-with-respect-to-R}). Hence by Theorem \ref{theor-Scott-functions-H-sober-1} (for $\mbox{H}=\mathsf{WD}$ or $\mbox{H}=\mathsf{RD}$), the Scott topology on $\mathbb{C}(X, Y)$ does not have property S with respect to $\mathsf{WD}$ and $\mathsf{RD}$.
\end{remark}

\vspace{3mm}

The following example gives an affirmative answer to Question C.

\begin{example}\label{exam-order-compatible-d-topology-on-complete-lattice-is-not-WF} (Johnstone's dcpo adding a top element)  Let $\mathbb{J}=\mathbb{N}^+\times (\mathbb{N}^+\cup \{\infty\})$ with ordering defined by $(j, k)\leq (m, n)$ if{}f $j = m$ and $k \leq n$, or $n =\infty$ and $k\leq m$ (see Figure 2). $\mathbb{J}$ is a well-known dcpo constructed by Johnstone in \cite{Johnstone-1981}.

\begin{figure}[ht]
	\centering
	\includegraphics[height=5cm,width=7cm]{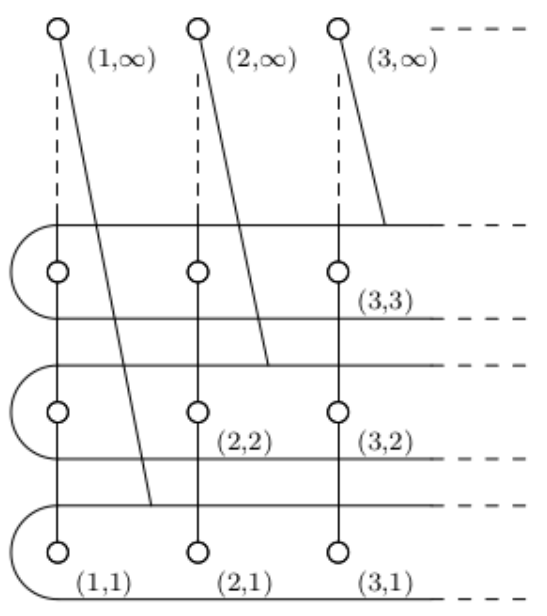}
	\caption{Johnstone's dcpo $\mathbb{J}$}
\end{figure}

The set $\mathbb{J}_{max}=\{(n, \infty) : n\in\mathbb{N}^+ \}$ consists of all maximal elements of $\mathbb{J}$. Adding a top element $\top$ to $\mathbb{J}$  yields a dcpo $\mathbb{J}_{\top}=\mathbb{J}\cup \{\top\}$ ($x< \top$ for any $x\in \mathbb{J}$). Then $\top$ is the largest element of $\mathbb{J}_{\top}$ and $\{\top\}\in \sigma (\mathbb{J}_{\top})$. The following four conclusions about $\Sigma~\!\!\mathbb{J}$ are known (see, for example, \cite[Exercise 8.3.9]{Goubault-2013}, \cite[Example 3.1]{Lu-Li-2017} and \cite[Example 26]{Xu-Wen-Xi-2023}):
\begin{enumerate}[\rm (1)]
\item If $U\in\sigma(\mathbb{J})\setminus \{\emptyset\}$, then there is $m\in\mathbb{N}^+$ such that $\{(n, \infty) : n\in \mathbb{N}^+, m\leq n\}\subseteq U$.
\item $\ir_c (\Sigma~\!\!\mathbb{J})=\{\overline{\{x\}}=\da_{\mathbb{J}} x : x\in \mathbb{J}\}\cup \{\mathbb{J}\}$.
\item $Q(\Sigma~\!\!\mathbb{J})=(2^{\mathbb{J}_{max}} \setminus \{\emptyset\})\bigcup \mathbf{Fin}~\!\mathbb{J}$.
\item $\Sigma~\!\!\mathbb{J}$ is not well-filtered.
\end{enumerate}
 Let $\tau=\{\emptyset\}\cup \{U\cup\{\top\} : U\in \sigma(\mathbb{J})\setminus \{\emptyset\}\}$. Then by (1) $\tau$ is a topology on  $\mathbb{J}_{\top}$ and $\tau=\sigma(\mathbb{J}_\top)\setminus \{\{\top\}\}$.
\begin{enumerate}[\rm (a)]
\item $\ir_c (\Sigma~\!\!\mathbb{J}_\top)=\{\overline{\{x\}}=\da_{\mathbb{J}_{\top}} x : x\in \mathbb{J}_\top\}\cup \{\mathbb{J}\}$ by (2).
\item $Q(\Sigma~\!\!\mathbb{J}_\top)=\{\ua_{\mathbb{J}_{\top}} G : G \mbox{~is nonempty and~} G\subseteq \mathbb{J}_{max}\cup\{\top\}\}\bigcup \mathbf{Fin}~\!\mathbb{J}_\top$ by (3).

\item $\upsilon(\mathbb{J}_{\top})\subsetneqq  \tau \subsetneqq \sigma(\mathbb{J}_{\top})$.

As $\tau=\sigma(\mathbb{J}_\top)\setminus \{\{\top\}\}$ and $\{\top\}\in \sigma(\mathbb{J}_\top)$, we have $\tau \subsetneqq  \sigma(\mathbb{J}_{\top})$. Clearly, $\{\top\}\not\in \upsilon(\mathbb{J}_{\top})$, whence $\upsilon(\mathbb{J}_{\top})\subseteq \sigma(\mathbb{J}_{\top})\setminus \{\{\top\}\}=\tau$. Let $A=\{(1, n) : n\in\mathbb{N}^+\}$. Then $A$ is a Scott closed subset of $\mathbb{J}_\top$, but it is not closed in $(\mathbb{N}^+, \upsilon(\mathbb{N}^+))$. Hence $\upsilon(\mathbb{J}_{\top})\subsetneqq  \tau$.

\item $(\mathbb{J}_{\top}, \tau)$ is sober.

By (2) and $\tau=\sigma(\mathbb{J}_\top)\setminus \{\{\top\}\}$, $\ir_c((\mathbb{J}_{\top}, \tau))= \bigl\{\mbox{cl}_{\sigma(\mathbb{J}_\top)}{\{x\}} : x \in \mathbb{J}\bigr\} \cup \{\mathbb{J}_{\top}\}=\bigl\{\mbox{cl}_{\sigma(\mathbb{J}_\top)}{\{x\}} : x \in \mathbb{J}_\top\}$. Thus $(\mathbb{J}_{\top}, \tau)$ is sober.

\item $\Sigma~\!\!\mathbb{J}_\top$ is not well-filtered.

Indeed, let $\mathcal K=\{\ua_{\mathbb{J}_\top} (\mathbb{J}_{max}\setminus F) : F\in (\mathbb{J}_{max})^{(<\omega)}\}$. Then by (b), $\mathcal K\subseteq Q(\Sigma~\!\!\mathbb{J}_\top)$ is a filtered family and $\bigcap\mathcal{K}=\bigcap_{F\in (\mathbb{J}_{max})^{(<\omega)}} \ua_{\mathbb{J}_\top} (\mathbb{J}_{max}\setminus F)=\bigcap_{F\in (\mathbb{J}_{max})^{(<\omega)}} ((\mathbb{J}_{max}\setminus F)\cup\{\top\})=\{\top\}\cup(\mathbb{J}_{max}\setminus \bigcup (\mathbb{J}_{max})^{(<\omega)})=\{\top\}\in \sigma (\mathbb{J}_\top)$, but there is no $F\in (\mathbb{J}_{max})^{(<\omega)}$ with $\ua_{\mathbb{J}_\top} (\mathbb{J}_{max}\setminus F)\subseteq \{\top\}$. Therefore, $\Sigma~\!\!\mathbb{J}_\top$ is not well-filtered.

\end{enumerate}
\end{example}

By Remark \ref{rem-problem-3-and-problem-2} and Example \ref{exam-order-compatible-d-topology-on-complete-lattice-is-not-WF}, we get a negative answer to Question B. More precisely, there are a sober space $X$ and a well-filtered space $Y$ such that the function space $[X\rightarrow Y]_\Sigma$ is not well-filtered, or equivalently, the Scott topology on $\mathbb{C}(X, Y)$ does not have property S with respect to $\mathsf{WD}$ and $\mathsf{RD}$.

\end{document}